\documentclass{amsart}

\usepackage[utf8]{inputenc} 
\usepackage[english]{babel} 
\usepackage{mathrsfs}

\theoremstyle{plain} 
\newtheorem{thm}{Theorem} 
\newtheorem{lem}[thm]{Lemma} 
\newtheorem{cor}[thm]{Corollary}

\theoremstyle{definition}

\theoremstyle{remark} 
\newtheorem*{rem}{Remark}

\begin{document} 
\title{Zeros of functions in Bergman--type Hilbert Spaces of Dirichlet Series} 
\author{Ole Fredrik Brevig} 
\date{\today} \address{Department of Mathematical Sciences, Norwegian University of Science and Technology (NTNU), NO-7491 Trondheim, Norway} \email{ole.brevig@math.ntnu.no} 
\begin{abstract}
	For a real number $\alpha$ the Hilbert space $\mathscr{D}_\alpha$ consists of those Dirichlet series $\sum_{n=1}^\infty a_n/n^s$ for which $\sum_{n=1}^\infty |a_n|^2/[d(n)]^\alpha < \infty$, where $d(n)$ denotes the number of divisors of $n$. We extend a theorem of Seip on the bounded zero sequences of functions in $\mathscr{D}_\alpha$ to the case $\alpha>0$. Generalizations to other weighted spaces of Dirichlet series are also discussed, as are partial results on the zeros of functions in the Hardy spaces of Dirichlet series $\mathscr{H}^p$, for $1\leq p <2$. 
\end{abstract}
\maketitle

\section{Introduction} Let $d(n)$ denote the divisor function let $\alpha$ be a real number. We are interested in the following Hilbert spaces of Dirichlet series:
\[\mathscr{D}_{\alpha} = \left\{f(s) = \sum_{n=1}^\infty \frac{a_n}{n^s} \, : \, \|f\|_{\mathscr{D}_{\alpha}}^2 = \sum_{n=1}^\infty \frac{|a_n|^2} {[d(n)]^\alpha}<\infty\right\}.\]
The functions of $\mathscr{D}_\alpha$ are analytic in $\mathbb{C}_{1/2} = \{s = \sigma+it \, : \, \sigma > 1/2\}$. Bounded Dirichlet series are almost periodic, and this implies that they have either no zeros or infinitely many zeros, as observed by Olsen and Seip in \cite{olsenseip}. This leads us to restrict our investigations to bounded zero sequences for spaces of Dirichlet series. In \cite{seipzero}, Seip studied bounded zero sequences for $\mathscr{D}_\alpha$, when $\alpha\leq 0$. This includes the Hardy--type $(\alpha=0)$ and Dirichlet--type $(\alpha<0)$ spaces. The topic of the present work is the Bergman--type spaces $(\alpha>0)$.

Let us therefore introduce the weighted Bergman spaces in the half-plane, $A_\beta$. For $\beta > 0$, these spaces consists of functions $F$ which are analytic in $\mathbb{C}_{1/2}$ and satisfy
\[\|F\|_{A_\beta} = \left(\int_{\mathbb{C}_{1/2}} |F(s)|^2 \left(\sigma-\frac{1}{2}\right)^{\beta-1}\, dm(s)\right)^\frac{1}{2}<\infty.\]
It was shown by Olsen in \cite{olsenlocal} that the local behavior of the spaces $\mathscr{D}_\alpha$ are similar to the spaces $A_\beta$, where $\beta=2^{\alpha}-1$. This relationship between $\alpha$ and $\beta$ will be retained throughout this paper.

For a class of analytic functions $\mathscr{C}$ on some domain $\Omega\subseteq \mathbb{C}$, we will say that a sequence $S$ of not necessarily distinct numbers in $\Omega$ is a zero sequence for $\mathscr{C}$ if there is some non-trivial $F \in \mathscr{C}$ vanishing on $S$, taking into account multiplicities. We will let $Z(\mathscr{C})$ denote the set of all zero sequences for $\mathscr{C}$.

A result proved by Horowitz in \cite{horowitz} shows that if $\mathscr{C}=A_\beta$ we may assume that $F$ vanishes precisely on $S\in Z(A_\beta)$, i.e. $F$ has no extraneous zeros in $\mathbb{C}_{1/2}$. We will exploit this fact to prove our main result. 
\begin{thm}
	\label{thm:zset} Suppose $S = (\sigma_j+it_j)$ is a bounded sequence of points in $\mathbb{C}_{1/2}$ and that $\alpha > 0$. Then there is a non-trivial function in $\mathscr{D}_{\alpha}$ vanishing on $S$ if and only if $S \in Z(A_\beta)$. 
\end{thm}
The ``only if'' part follows from the local embedding of $\mathscr{D}_\alpha$ into $A_\beta$ of Theorem 1 and Example 4 from \cite{olsenlocal}. To prove the ``if'' part, we will adapt the methods of \cite{seipzero}, where an analogous result for $\alpha \leq 0$ was obtained. 

The ``if'' part can essentially be split into two steps. The first step is a discretization lemma, which depends on the properties of $\mathscr{D}_\alpha$ --- or rather the weights $[d(n)]^\alpha$. The second step is an iterative scheme, where the properties of $A_\beta$ become more prominent. 

Comparing this with \cite{seipzero}, the first step is somewhat harder, since we require very precise estimates on the weights as $\alpha$ grows to infinity. The second step is considerably easier, mainly due to the fact that the norms of $A_\beta$ are easier to work with than those of the Dirichlet spaces used in \cite{seipzero}.

We will use the notation $f(x) \ll g(x)$ to indicate that there is some constant $C>0$ so that $|f(x)| \leq C g(x)$. Sometimes the constant $C$ may depend on certain parameters, and this will be specified in the text. Moreover, we write $f(x) \asymp g(x)$ if both $f(x) \ll g(x)$ and $g(x) \ll f(x)$ hold.

\section{Proof of Theorem \ref{thm:zset}}
We begin with the Paley--Wiener representation of functions $F \in A_\beta$, and seek to construct a Dirichlet series $f \in \mathscr{D}_\alpha$ which approximates $F$. 
\begin{lem}
	[Paley--Wiener Representation] \label{lem:paleywiener} $A_\beta$ is isometrically isomorphic to
	\[L_\beta^2 = \left\{\phi \text{ measurable on } [0,\infty) \, : \, \|\phi\|_{L^2_\beta}^2 = \frac{2\pi\Gamma(\beta)}{2^\beta} \int_0^\infty |\phi(\xi)|^2 \,\frac{d\xi}{\xi^\beta}<\infty\right\},\]
	under the Laplace transformation
	\[F(s) = \int_0^\infty \phi(\xi) e^{-(s-1/2)\xi}\, d\xi.\]
	\begin{proof}
		A proof can be found in \cite{paleywiener}. 
	\end{proof}
\end{lem}
The other ingredient needed for the discretization lemma is estimates on the growth of $[d(n)]^\alpha$. We will partition the integers into blocks and use an average order type estimate. To prove this estimate, we will need the precise form of a formula stated by Ramanujan \cite{ramanujan} and proved by Wilson \cite{wilson}: For any real number $\alpha$ and any integer $\nu>2^{\alpha}-2$, we have 
\begin{equation}
	D_\alpha(x) = \sum_{n\leq x} [d(n)]^\alpha = x (\log{x})^{2^\alpha-1} \left(\sum_{\lambda=0}^\nu \frac{A_\lambda}{(\log{x})^\lambda}+\mathcal{O}\left(\frac{1}{(\log{x})^{\nu+1}}\right)\right). \label{eq:ramaD} 
	\end{equation}
Wilson's proof of \eqref{eq:ramaD} can be considered at special case of Selberg--Delange method. For more about the Selberg--Delange method, we refer to Chapter II.5 of \cite{tenenbaum}. However, we mention that the coefficients $A_\lambda$ depend on the coefficients of the Dirichlet series $\phi_\alpha$, which we implicitly define through the relation
\begin{equation}
	\zeta_\alpha(s) = \sum_{n=1}^\infty [d(n)]^\alpha n^{-s} = \prod_{j=1}^\infty \left(1 + \sum_{k=1}^\infty (k+1)^\alpha p_j^{-sk} \right) = \left[\zeta(s)\right]^{2^\alpha} \phi_\alpha(s).
	\label{eq:phidef}
\end{equation}
The partial sums of the coefficients of $\zeta_\alpha$ are estimated through Perron's formula and the residue theorem. While \eqref{eq:phidef} is only valid for $\Re(s)>1$, a simple computation using Euler products shows that $\phi_\alpha$ converges for $\Re(s)>1/2$, and thus Theorem 5 of \cite{tenenbaum} may be applied. In particular, the coefficients $A_\lambda$ depend on the coefficients of $\phi_\alpha$, and since the coefficients of $\phi_\alpha$ depend continuously on $\alpha$, so does $A_\lambda$ in \eqref{eq:ramaD}. 
\begin{lem}
	\label{lem:dirisum} Let $\alpha$ be a real number and $0 < \gamma < 1$. Then 
	\begin{equation}
		\label{eq:divsum} \sum_{j^\gamma \leq \log{n} \leq (j+1)^\gamma} \frac{[d(n)]^\alpha}{n} \asymp j^{\gamma2^\alpha-1}, 
	\end{equation}
	as $j \to \infty$. The implied constants may depend on $\alpha$ and $\gamma$. 
	\begin{proof}
		We will first assume that $2^\alpha$ is not an integer. Fix $\nu$ such that $\nu > 2^{\alpha}-1$ and $\nu >1/\gamma-1$. We use Abel summation to rewrite 
		\begin{equation}
			\label{eq:abel} \sum_{y < n \leq x} \frac{[d(n)]^\alpha}{n} = \frac{D_\alpha(x)}{x}-\frac{D_\alpha(y)}{y}+ \int_y^x \frac{D_\alpha(z)}{z^2}\, dz. 
		\end{equation}
		By using \eqref{eq:ramaD} and the fact that $2^{\alpha}-1-\nu<0$ we perform some standard calculations to estimate 
		\begin{align*}
			\frac{D_\alpha(x)}{x}-\frac{D_\alpha(y)}{y} &= \sum_{\lambda=0}^\nu A_\lambda \left((\log{x})^{2^\alpha-1-\lambda}-(\log{y})^{2^\alpha-1-\lambda}\right) + \mathcal{O}\left((\log{y})^{2^\alpha-2-\nu}\right), \\
			\int_y^x \frac{D_\alpha(z)}{z^2}\,dz &= \sum_{\lambda=0}^\nu \frac{A_\lambda}{2^\alpha-\lambda}\left((\log{x})^{2^\alpha-\lambda}-(\log{y})^{2^\alpha-\lambda}\right) + \mathcal{O}\left((\log{y})^{2^\alpha-1-\nu}\right). 
		\end{align*}
		Let us now take $x=\exp\left((j+1)^\gamma\right)$ and $y = \exp\left(j^{\gamma}\right)$. For any exponent $\eta$ it is clear that
		\[(\log{x})^\eta-(\log{y})^\eta = \gamma \eta j^{\gamma\eta-1}\left(1+\mathcal{O}\left(\frac{1}{j}\right)\right).\]
		Hence we have 
		\begin{align*}
			\frac{D_\alpha(x)}{x}-\frac{D_\alpha(y)}{y} &\asymp \sum_{\lambda=0}^\nu A_\lambda(\gamma(2^{\alpha}-1-\lambda)) j^{\gamma(2^\alpha-1-\lambda)-1} + \mathcal{O}\left(j^{\gamma(2^\alpha-2-\nu)}\right), \\
			\int_y^x \frac{D_\alpha(z)}{z^2}\,dz &\asymp \sum_{\lambda=0}^\nu A_\lambda j^{\gamma(2^\alpha-\lambda)-1}+ \mathcal{O}\left(j^{\gamma(2^\alpha-1-\nu)}\right). 
		\end{align*}
		We combine these estimates with \eqref{eq:abel} to obtain
		\begin{equation}
			\sum_{j^\gamma \leq \log{n} \leq (j+1)^\gamma} \frac{[d(n)]^\alpha}{n} \asymp j^{\gamma2^\alpha-1} \left( A_0 + \sum_{\lambda=1}^\nu \frac{B_\lambda}{j^{\gamma\lambda}}+\mathcal{O}\left(\frac{1}{j^{\gamma2^\alpha-1-\gamma(2^\alpha-1-\nu)}}\right)\right),
			\label{eq:finest}
		\end{equation}
		where $B_\lambda = A_\lambda + A_{\lambda-1} \gamma \left(2^\alpha-\lambda\right)$. This proves \eqref{eq:divsum} since $\nu>1/\gamma-1$. By continuity on both sides of \eqref{eq:finest}, the assumption that $2^\alpha$ is not an integer may be dropped.
	\end{proof}
\end{lem}
The parameter $0 < \gamma < 1$ will be used to control the ``block size'' in our partition of the integers. It will become apparent that as $\alpha$ grows to infinity, we must be able to let $\gamma$ tend to $0$. In \cite{seipzero} it was sufficient to have a similar estimate only for $1/2 < \gamma < 1$. 
\begin{lem}[Discretization Lemma] \label{lem:disc} Let $\alpha>0$ and let $N$ be a sufficiently large positive integer. Then there exists positive constants $A$ and $B$ (depending on $\alpha$, but not $N$) such that the following holds: For every function $\phi \in L_\beta^2$ supported on $\left[\log{N},\infty\right)$, there is a function of the form
	\[f(s) = \sum_{n=N}^\infty \frac{a_n}{n^s}\]
	in $\mathscr{D}_\alpha$ such that $\|f\|_{\mathscr{D}_\alpha} \leq A \|\phi\|_{L^2_\beta}$. Moreover, $f$ may be chosen so that
	\[\Phi(s) = \int_{\log{N}}^\infty \phi(\xi) e^{-(s-1/2)\xi}\,d\xi - f(s)\]
	enjoys the estimate 
	\[|\Phi(s)| \leq B |s-1/2| N^{-\sigma+1/2}(\log{N})^{-1}\|\phi\|_{L^2_\beta},\]
	in $\mathbb{C}_{1/2}$.
	\begin{proof}
		Let $\gamma = 2/(4+2^\alpha)$ and let $J$ be the largest integer smaller than $(\log(N))^{1/\gamma}$. For $j\geq J$, let $n_j$ be the smallest integer $n$ such that $e^{j^\gamma}\leq n$. When $\gamma$ is small it is possible that $n_j = n_{j+1}$. This can be avoided by taking $N$ sufficiently large. Set $\xi_{n_j}=j^\gamma$ and for $n_j < n \leq n_{j+1}$ iteratively choose $\xi_n$ such that 
		\begin{equation}
			\label{eq:iteprod} \frac{\xi_{n+1}^{\beta+1}-\xi_n^{\beta+1}}{\beta+1} = A_j \frac{[d(n)]^\alpha}{n}, 
		\end{equation}
		where $A_j$ is chosen so that $\xi_{n_{j+1}}=(j+1)^\gamma$. Clearly, Lemma \ref{eq:divsum} implies that $A_j$ is bounded as $j \to \infty$. Let us set
		\[a_n = \sqrt{n} \int_{\xi_n}^{\xi_{n+1}} \phi(\xi)\, d\xi.\]
		A simple computation using the Cauchy--Schwarz inequality shows that
		\[|a_n|^2 = n \left|\int_{\xi_n}^{\xi_{n+1}}\phi(\xi) \,d\xi\right|^2 \leq n \cdot\frac{\xi_{n+1}^{\beta+1}-\xi_n^{\beta+1}}{\beta+1} \int_{\xi_n}^{\xi_{n+1}}|\phi(\xi)|^2\, \frac{d\xi}{\xi^\beta}.\]
		In view of \eqref{eq:iteprod} it is clear that $\|f\|_{\mathscr{D}_{\alpha}}\leq A \|\phi\|_{L_\beta^2}$. Now, if $n_j \leq n \leq n_{j+1}$ and $\xi \in [\xi_{n_j},\, \xi_{n_{j+1}}]$ we see that 
		\begin{equation}
			\label{eq:intest} \left|e^{-(s-1/2)}-n^{-(s-1/2)}\right| \leq N^{-\sigma+1/2} |s-1/2| j^{\gamma-1}. 
		\end{equation}
		Then, by \eqref{eq:intest} and the Cauchy--Schwarz inequality 
		\begin{align*}
			&|\Phi(s)| \leq N^{-\sigma+1/2}|s-1/2|\sum_{j=J}^\infty j^{\gamma-1}\sum_{n = n_j}^{n_{j+1}-1} \left(\frac{\xi_{n+1}^\beta- \xi_{n}^\beta}{\beta}\right)^\frac{1}{2}\left(\int_{\xi_n}^{\xi_{n+1}} |\phi(\xi)|^2 \, \frac{d\xi}{\xi^\beta}\right)^\frac{1}{2}. \intertext{By using the Cauchy--Schwarz inequality again with \eqref{eq:iteprod} we get} &|\Phi(s)| \ll N^{-\sigma+1/2}|s-1/2| \sum_{j=J}^\infty j^{\gamma-1} \left(\sum_{n=n_j}^{n_{j+1}-1} \frac{[d(n)]^\alpha}{n}\right)^\frac{1}{2}\left(\int_{\xi_{n_j}}^{\xi_{n_{j+1}}}|\phi(\xi)|^2\, \frac{d\xi}{\xi^\beta}\right)^\frac{1}{2}. \intertext{Now Lemma \ref{lem:dirisum} and the Cauchy--Schwarz inequality yield}  &|\Phi(s)| \ll N^{-\sigma+1/2} |s-1/2| \left(\sum_{j=J}^\infty j^{(2+2^\alpha)\gamma-3}\right)^\frac{1}{2} \Bigg(\int_{\log{N}}^\infty |\phi(\xi)|^2\, \frac{d\xi}{\xi^\beta}\Bigg)^\frac{1}{2}. 
		\end{align*}
		The series converges since $\gamma < 2/(2+2^\alpha)$. The proof is completed by a standard estimate of the convergent series,
		\[\left(\sum_{j=J}^\infty j^{(2+2^\alpha)\gamma-3}\right)^\frac{1}{2} \ll (\log{N})^{((2+2^\alpha)\gamma-2)/(2\gamma)} = (\log{N})^{-1},\]
		where we used that $J \asymp (\log{N})^{1/\gamma}$. 
	\end{proof}
\end{lem}
The final result needed for the iterative scheme is the following simple lemma on the $\overline{\partial}$-equation. We omit the proof, which is obvious. 
\begin{lem}
	\label{lem:dbar} Suppose $g$ is a continuous function on $\mathbb{C}_{1/2}$, supported on
	\[\Omega(R,\tau) = \{s=\sigma + i t\, : \, 1/2 \leq \sigma \leq 1/2+\tau,\, -R \leq t \leq R\},\]
	for some positive real numbers $\tau$ and $R$. Then
	\[u(s) = \frac{1}{\pi} \int_\Omega \frac{g(w)}{s-w}\, dm(w)\]
	solves $\overline{\partial}u = g$ in $\mathbb{C}_{1/2}$ and satisfies $\|u\|_\infty \leq C_\Omega \|g\|_\infty$.
\end{lem}

We have now collected all our preliminary results and are ready to begin the proof of Theorem \ref{thm:zset}. For any positive integer $N$ we set $E_N(s) = N^{-s+1/2}$ and consider the space $E_N A_\beta$. By a substitution it is evident that any $F \in E_N A_\beta$ can be represented as
\[F(s) = \int_{\log{N}}^\infty \phi(\xi)e^{-(s-1/2)\xi}\, d\xi\]
for some $\phi \in L_\beta^2[\log{N},\infty)$, in view of Lemma \ref{lem:paleywiener}. 
\begin{proof}
	[Final step in the proof of Theorem \ref{thm:zset}] Let us fix $\alpha>0$ and a bounded sequence $S=(\sigma_j+it_j)\in Z(A_\beta)$. From this point all constants may depend on $\alpha$ and $S$. Since $S$ is bounded we may assume $S \subset \Omega(R-2,\tau-2)$ for some $R, \, \tau > 2$. Let $\Theta$ be some smooth function defined on $\overline{\mathbb{C}_{1/2}}$ with the following properties: 
	\begin{itemize}
		\item $\Theta$ is supported on $\Omega(R,\tau)$, 
		\item $\Theta(s)=1$ for $s \in \Omega(R-1,\tau-1)$, 
		\item $|\overline{\partial}\Theta(s)|\leq 2$. 
	\end{itemize}
	Let $G \in A_\beta$ vanish precisely on $S$ and assume furthermore that $\|G\|_{A_\beta}=1$. Now, suppose that $F \in E_N A_\beta$, and let $f \in \mathscr{D}_\alpha$ be the function obtained by applying Lemma \ref{lem:disc} to $F$, and $\Phi = F - f$. Moreover, let $u$ denote the solution to the equation 
	\begin{equation}
		\label{eq:udb} \overline{\partial} u = \frac{\overline{\partial}(\Theta \Phi)}{G E_N}. 
	\end{equation}
	The right hand side of \eqref{eq:udb} is a smooth function compactly supported on $\Omega(R,\tau)$ since $|G(s)|$ is bounded from below where $\overline{\partial} \Theta(s)\neq 0$. We can use Lemma \ref{lem:dbar} and Lemma \ref{lem:paleywiener} to estimate 
	\begin{equation}
		\label{eq:uest} \|u\|_\infty \ll \left\|\frac{\overline{ 
		\partial}(\Theta \Phi)}{G E_N}\right\|_\infty \ll (\log{N})^{-1} \|\phi\|_{L_\beta^2} = (\log{N})^{-1} \|F\|_{A_\beta}. 
	\end{equation}
	We set $T_NF= \Theta \Phi-G E_N u$. The function $T_N F$ has the following properties: 
	\begin{itemize}
		\item $T_NF(s)= \Phi(s)$ for $s \in S$, 
		\item $T_NF$ is analytic in $\mathbb{C}_{1/2}$ since $\overline{ 
		\partial} T_NF(s)=0$ for $s \in \mathbb{C}_{1/2}$, 
		\item $T_NF \in E_N A_\beta$, by the compact support of $\Theta$ and the estimate \eqref{eq:uest}. 
	\end{itemize}
	Hence $T_N$ defines an operator on $E_N A_\beta$. By the triangle inequality, Lemma \ref{lem:disc} and the fact that $\Theta$ has compact support, it is clear that
	\[\|T_N F \|_{A_\beta} \leq \|\Theta \Phi \|_{A_\beta} + \|GE_N u\|_{A_\beta} \ll (\log{N})^{-1} \|\phi\|_{L_\beta^2} + \|u\|_\infty\|G\|_{A_\beta}.\]
	Since $\|G\|_{A_\beta}=1$ and $\|\phi\|_{L_\beta^2}=\|F\|_{A_\beta}$ we have $\|T_N\| \ll (\log{N})^{-1}$ in view of \eqref{eq:uest}. Let $N$ be large, but arbitrary, and define $F_0(s)= E_N(s)G(s)$. Then $F_0 \in E_N A_\beta$ and its norm in this space is $\leq 1$. Set
	\[F_j = T_N^j F_0.\]
	Let $f_j$ be the Dirichlet series of Lemma \ref{lem:disc} obtained from $F_j$. Then $f_0+F_1$ vanishes on $S$, since
	\[f_0(s) + F_1(s) = f_0(s) + T_N F_0(s) = f_0(s) + F_0(s) - f_0(s) = F_0(s) = 0,\]
	for $s\in S$, by the fact that $T_N F(s) = \Phi(s)$ for $s \in S$. Iteratively, the function $f_0+f_1+\cdots+f_j + F_{j+1}$ also vanishes on $S$. Define
	\[f(s) = \sum_{j=0}^\infty f_j(s)\]
	and choose $N$ so large that $\|T_N\|<1$ so that $\|F_j\|_{A_\beta} \to 0$ and, say
	\[|f(1)| > \sum_{j=1}^\infty |f_j(1)|,\]
	so that $f$ is non-trivial in $\mathscr{D}_\alpha$ and vanishing on $S$. 
\end{proof}

By again following \cite{seipzero}, we can modify the iterative scheme in the following way: Let $F \in A_\beta$ be arbitrary, and set $F_0 = F$. Using the algorithm in the same manner as above, we see that $F_1(s) + f_0(s) = F_0(s)$ for $s \in S$. Moreover,
\[F_{j+1}(s) + f_j(s) + f_{j-1}(s) + \cdots + f_0(s) = F(s),\]
for $s \in S$. Continuing as above, we obtain the following result: 
\begin{cor} \label{cor:interpolation} 
	Suppose $S = (\sigma_j+it_j) \in Z(A_\beta)$ is bounded. For every function $F \in A_\beta$ there is some $f\in \mathscr{D}_{\alpha}$ such that $f(s)=F(s)$ on $S$. 
\end{cor}

We can extend Theorem \ref{thm:zset} and Corollary \ref{cor:interpolation} by considering different weights. Let $w=(w_1,\,w_2,\,\ldots\,)$ be a non-negative weight. Define the Hilbert space of Dirichlet series $\mathscr{D}_w$ in the same manner as above, with the added convention that the basis vector $n^{-s}$ is excluded if $w_n=0$. Theorem 1 in \cite{olsenlocal} states that $\mathscr{D}_w$ embeds locally into $A_\beta$ if and only if 
\begin{equation}
	\label{eq:embedest} \sum_{n\leq x} w_n \ll x(\log{x})^{\beta}, 
\end{equation}
where $\beta>0$. By modifying the proof of our Theorem \ref{thm:zset}, we can obtain a similar result for $\mathscr{D}_w$ with respect to $A_\beta$ provided we additionally have 
\begin{equation}
	\label{eq:ourest} \sum_{j^\gamma \leq \log{n} \leq (j+1)^\gamma} \frac{w_n}{n} \asymp j^{\gamma(\beta+1)-1}, 
\end{equation}
as $j \to \infty$, for some $0 < \gamma < 2/(3+\beta)$. Several of the weights considered in \cite{olsenlocal} are possible, but we only mention the case $w_n = (\log{n})^\beta$ for $\beta>0$. These spaces were introduced by McCarthy in \cite{mccarthy}. It is easy to show that these weights satisfy \eqref{eq:embedest} and \eqref{eq:ourest} for any $0 < \gamma < 1$, and similar results with respect to $A_\beta$ are obtained.
\begin{rem}
	The embeddings of \cite{olsenlocal} extend to any $\beta \leq 0$, in view of \eqref{eq:embedest}, and we get the Hardy space $(\beta=0)$ and Dirichlet spaces $(\beta<0$) in the half-plane. We can extend the results in \cite{seipzero} in a similar manner as above. However, this is only possible for $-1 \leq \beta < 0$. The method of \cite{seipzero} breaks down for $\beta < -1$ due to the fact that the norms of the corresponding Dirichlet spaces in the half-plane uses higher order derivatives and different estimates are needed.
\end{rem}

\section{Blaschke-type conditions for $\mathscr{D}_\alpha$ and $\mathscr{H}^p$}
Now that we have identified the bounded zero sequences of $\mathscr{D}_\alpha$ as those of $A_\beta$, let us consider necessary and sufficient conditions for bounded zero sequences of $A_\beta$. The zero sequences of Bergman spaces in the unit disc $\mathbb{D}$ have attracted considerable attention. We refer to the monograph \cite{HKZ}. For $\beta>0$, these are the spaces
\[A_\beta(\mathbb{D}) = \left\{F \in H(\mathbb{D}) \,:\, \|F\| = \int_{\mathbb{D}} |F(z)|^2(1-|z|)^{\beta-1}dm(z) < \infty\right\}.\]
Results pertaining to zero sequences of $A_\beta(\mathbb{D})$ are relevant to our case since 
\[\phi(s)=\frac{s-3/2}{s+1/2}\]
is a conformal mapping from $\mathbb{C}_{1/2}$ to $\mathbb{D}$, and
\[F \mapsto (s+1/2)^{-2(\beta+1)} F\left(\frac{s-3/2}{s+1/2}\right)\]
defines an isometric isomorphism from $A_\beta(\mathbb{D})$ to $A_\beta$. This implies that $S \in Z(A_\beta)$ if and only if $\phi(S) \in Z(A_\beta(\mathbb{D}))$. Since the Hardy space $H^2(\mathbb{D})$ is included in $A_\beta(\mathbb{D})$ for every $\beta>0$, it is clear that the Blaschke condition
\begin{equation}
	\label{eq:blaschke} \sum_{j} (\sigma_j-1/2) < \infty 
\end{equation}
is sufficient for bounded zero sequences of $A_\beta$. Moreover, Theorem 4.1 of \cite{HKZ} shows that the Blaschke condition \eqref{eq:blaschke} is both necessary and sufficient provided the bounded sequence $S$ is contained in any cone $|t-t_0|\leq c(\sigma-1/2)$. Unfortunately, the situation becomes more complicated in the general case and we do not have a precise Blaschke-type condition for bounded zero sequences. In fact, for every $\epsilon>0$ and every $A_\beta$ a necessary condition for bounded zero sequences is
\begin{equation} \label{eq:blepsilon}
	\sum_{j} (\sigma_j-1/2)^{1+\epsilon}<\infty,
\end{equation}
by Corollary 4.8 of \cite{HKZ}. Clearly, this condition does not offer any insight into what happens as $\beta \to 0^+$. However, using the notion of density introduced by Korenblum in \cite{korenblum} it is possible to provide a generalized condition describing the geometrical information of the zero sequences of $A_\beta(\mathbb{D})$. The most precise results on Korenblum's density are obtained by Seip in \cite{seip1995korenblum}. We omit the details, only mentioning that this generalized condition in a certain sense tends to \eqref{eq:blaschke} when $\beta \to 0^+$.

The Hardy spaces of Dirichlet series $\mathscr{H}^p$, $1 \leq p < \infty$, can be defined as the closure of the set of all Dirichlet polynomials with respect to the norms
\[\left\|\sum_{n=1}^N \frac{a_n}{n^s}\right\|_{\mathscr{H}^p} = \lim_{T \to \infty}\left(\frac{1}{2T} \int_{-T}^T \left|\sum_{n=1}^N \frac{a_n}{n^{it}}\right|^p\,dt\right)^\frac{1}{p}.\]
For the basic properties of these spaces we refer to \cite{HLS} and \cite{bayart}. However, we immediately observe that $\mathscr{H}^2 = \mathscr{D}_0$. In \cite{seipzero}, the bounded zero sequences of the spaces $\mathscr{H}^p$, for $2\leq p < \infty$, are studied. In particular, for $\mathscr{H}^2$ the Blaschke condition \eqref{eq:blaschke} is shown to be both necessary and sufficient. Results for $2<p<\infty$ are obtained through embeddings $\mathscr{D}_\alpha \subset \mathscr{H}^p \subset \mathscr{H}^2$, where $\alpha<0$ depends on $p$. The embedding of $\mathscr{H}^p$ into $\mathscr{H}^2$ implies that the Blaschke condition \eqref{eq:blaschke} is necessary for $\mathscr{H}^p$. 

The sufficient conditions are obtained through a similar result as Theorem \ref{thm:zset}: For $\alpha<0$, the spaces $\mathscr{D}_\alpha$ have the same bounded zero sequences as certain weighted Dirichlet spaces in $\mathbb{C}_{1/2}$. In particular, for $2 < p < \infty$ there is some $0 < \gamma < 1$ such that a sufficient condition for bounded zero sequences of $\mathscr{H}^p$ is
\begin{equation} \label{eq:diriblas}
	\sum_{j} (\sigma_j-1/2)^{1-\gamma}<\infty,
\end{equation}
and moreover $\gamma \to 0$ as $p \to 2^-$. We omit the details, which can be found in \cite{seipzero}.

We will now consider the case $1 \leq p < 2$. That $\mathscr{H}^2 \subset \mathscr{H}^p \subseteq \mathscr{H}^1$ for $1\leq p < 2$ is trivial, and this shows that \eqref{eq:blaschke} is a sufficient condition for bounded zero sequences of $\mathscr{H}^p$. In \cite{helson}, Helson proved the beautiful inequality 
\begin{equation} \label{eq:helson} 
	\|f\|_{\mathscr{D}_1} = \left(\sum_{n=1}^\infty \frac{|a_n|^2}{d(n)}\right)^\frac{1}{2} \leq \|f\|_{\mathscr{H}^1}, 
\end{equation}
which implies that $\mathscr{H}^p \subset \mathscr{D}_1$. This shows that the Blaschke-type condition \eqref{eq:blepsilon} is necessary for bounded zero sequences of $\mathscr{H}^p$, for every $\epsilon>0$. Regrettably, this means we are unable to specify how the situation changes as $p \to 2^-$, in a manner similar to \eqref{eq:diriblas}. However, if we again restrict $S$ to the cone $|t-t_0|\leq c(\sigma-1/2)$, the Blaschke condition \eqref{eq:blaschke} is both necessary and sufficient for bounded zero sequences of $\mathscr{H}^p$. 

\begin{rem}
	The Blaschke condition \eqref{eq:blaschke} is well-known to be necessary and sufficient for bounded zero sequences of the Hardy spaces $H^p(\mathbb{C}_{1/2})$. By a theorem in \cite{HLS}, $\mathscr{H}^2$ embeds locally into $H^2(\mathbb{C}_{1/2})$. This trivially extends to even integers $p$. Whether the local embedding extends to every $p \geq 1$ is an open question. Observe that if \eqref{eq:blaschke} is not the optimal necessary condition for bounded zero sequences of $\mathscr{H}^p$, when $1\leq p < 2$, then the local embedding would be impossible for these $p$. However, since \eqref{eq:diriblas} is a sufficient condition for bounded zero sequences of $\mathscr{H}^p$ when $p\geq 2$, its optimality would not contradict the local embedding for these $p$.
\end{rem}

\section*{Acknowledgements}
This paper constitutes a part of the author's PhD studies under the advice of Kristian Seip, whose feedback the author is grateful for. The author would also like to extend his gratitude to Jan-Fredrik Olsen for helpful discussions pertaining to Section 3.
\bibliographystyle{amsplain}
\bibliography{z}
\end{document}